
\baselineskip=14pt
\parskip=10pt

\magnification=\magstephalf

\def\N{{\cal N}}
\def\1{{\overline{1}}}
\def\2{{\overline{2}}}
\parindent=0pt
\overfullrule=0in

\def\frac#1#2{{#1 \over #2}}
\centerline
{\bf 
An Experimental (yet fully rigorous!) Study of a certain ``Measure Of Disarray" 
}
\centerline
{\bf that 12-year Noga Alon Proved was always Even}
\bigskip
\centerline
{\it Shalosh B. EKHAD and Doron ZEILBERGER}

\bigskip

{\bf Abstract}:  We study in depth a certain permutation statistic that was the subject of a brilliant insight by 12-year-old Noga Alon.
Our approach is purely empirical and experimental, yet it is fully rigorous, thereby debunking, yet another time, the myth that mathematics is always a deductive science.

{\bf Maple package}

This article is accompanied by a Maple package, {\tt Noga12.txt}, available free of charge from

{\tt https://sites.math.rutgers.edu/\~{}zeilberg/tokhniot/Noga12.txt} \quad.

The front of this article 

{\tt https://sites.math.rutgers.edu/\~{}zeilberg/mamarim/mamarimhtml/noga12.html} \quad ,

contains numerous output files that will be referred to later on.

{\bf Preface: How it all started}

In a beautiful new ``coffee table book", ``Do not Erase" [W], by the very talented artistic photographer Jessica Wynne, there are pictures of more than one hundred blackboards by 
a very diverse set of mathematicians. One of them is Noga Alon's blackboard. Each blackboard photo is accompanied by a short essay by the creator of the blackboard, where they often describe how 
they decided to become mathematicians. According to Noga Alon, the {\it epiphany}  occurred when he was twelve years old. Here are his exact words ([W], p. 52):

{\it \qquad \qquad ``One aspect of mathematics  that always fascinated me is its objectivity. When I was twelve years old, my parents had friends during the Eurovision song contest. They amused themselves
by trying to guess the final results of the contest, which are determined by the votes of the participating countries - they agreed on a formula (which I am not including here) for computing the  score for each guest
given the end results, and when the results were announced, they found out that all the scores were even integers. This led to a discussion about whether this must always be the case.
My mother (who had excellent mathematical intuition)  thought that it was indeed the case, but one of the guests, who was an engineer, insisted that this wasn't so. Since they were
unable to agree, my mother suggested asking me. Although I was quite young, she knew that I liked mathematics and believed that I would be able to settle the question. I thought about the
problem and realized that, indeed, the score must always be even.

The main reason that I remember the story, however, is not because I was able to find a proof, but because I managed to show the engineer that he was wrong. The fact that a twelve-year-old
can convince a grown engineer about such a statement is an impressive  demonstration of the objective nature of mathematics.''}

We were intrigued, and sent Noga Alon the following email message.

{\tt
I recently bought a beautiful book ``Do not Erase", by Jessica Wynne, and was
pleased to see that one of the blackboards is yours (btw, my Rutgers colleague, Jeff Kahn, should be honored
that he made it to your blackboard!). 

I also liked very much the accompanying essay where you mentioned
how, when you were twelve-year-old,  you managed to convince a grown-up engineer that he was wrong regarding
a certain  score related to the Eurovision contest.

What was the actual problem? What was your proof? 
}

Shortly after, Noga replied:

{\tt
 The details are as follows: Each
participant guessed a ranking of all countries, and the score of his/her
guess was the sum, over all countries, of the absolute value of the
difference between the guessed rank of each country and its final rank.
The question the guests of my parents had was whether or not the score must
always be even.

To show that indeed the score must be even note that the guessed ranking
is a permutation sigma of the numbers from 1 to n, (if there are n
countries).  The final ranking is also a permutation pi of these
numbers. The score is

(1) sum over i of | sigma(i)- pi (i)|

But if we forget the absolute values then clearly

sum over i of  (sigma (i)-pi (i))=0 (which is even)

and whenever we change a number (sigma (i)-pi (i)) to its
absolute value | sigma(i)- pi (i)| we change the sum by an even
number, showing that (1) is also even.

Not very deep, but as far as I remember I was quite happy noticing
it, especially since it has been simple enough and could be
explained to pretty much anybody who is willing to listen.
}

{\bf The Eurovision-Alon-Party `Measure of Disarray' was first invented in 1910}

This got us intrigued.
We first realized that it  is natural to make one of the permutations the identity permutation. In other words
rename the countries participating in the contest by their ranking, so the natural quantity, that we named $na(\pi)$ (in honor of Noga Alon), may be defined on a permutation
$\pi=[\pi_1, \dots, \pi_n]$ of  $\{1,\dots, n\}$  as follows:
$$
na(\pi):= \sum_{i=1}^{n} \, |\pi_i - i| \quad .
$$

Obviously this is a {\it measure} of `how scrambled' the permutation is. This is what our dear friend Dominique Foata calls a {\it permutation statistic}, a central object in
enumerative combinatorics.
But we never saw it before! Indeed the annoying {\bf absolute value} was dissonant to our mathematical ear. 

Nevertheless we started to investigate more closely this {\it measure of disarray}. 
In particular we found the first ten moments, as explicit expressions in $n$, and proved that the scaled moments (see below) tend, as $n$ goes to infinity,
 to those of the standard normal distribution $\N(0,1)$ (i.e  $1,0,3,0,15,0,105,0,945$), indicating
that it is, most probably, {\it asymptotically normal} (and we rigorously proved it for the first $10$ moments, also see below).

Before going public, we thought it prudent to email the great {\bf maven}, Persi Diaconis. On  Aug 26, 2021, 12:42 PM, one of us (DZ) sent him an email that included the following questions 
(regarding Noga Alon parents' permutation statistic):

{\tt

(i) Have this permutation statistic been considered before by statisticians?

 (ii) Is it known, or easily follows from general theorems, that it is asymptotically normal?

(Shalosh can prove that the scaled moments converge to those of N(0,1) up to
the tenth moment).
}

$25$ minutes later (Aug 26, 2021, 1:07 PM), we got the following reply:

{\tt
Hi Doron,

   Statisticians call this statistic `spearman's footrule'. Spearman was a psychologist who used it around 1910. And indeed I (personally) proved the limiting normality in 
a paper with Ron Graham (metrics on permutations).  I think that someone like Hausdorff (or Hurwitz) figured out the first two moments for Spearman. 
You might take a look at chapter 6 of my book which discusses other metrics and their limiting normality. It would be interesting to know the higher moments since one can use them to get 'Edgeworth corrections' to the CLT ...

Best wishes, Hope this helps. 

Persi Diaconis
}

We immediately looked up the paper that Persi Diaconis wrote with Ron Graham [DG], and found out that they used human ingenuity to figure out the expectation and variance.
They also applied a general deep theorem of
Hoeffding to prove {\it asymptotic normality}. 

In this paper we do the following two feats.

$\bullet$ Derive experimentally (yet rigorously) {\bf explicit} expressions for the first ten moments (about the mean), thereby giving a partial (but fully elementary) proof
of the asymptotic normality (proving that the scaled moments of Spearman's footrule converge to those of the standard normal distribution, 
$\frac{1}{\sqrt{2\pi}} \, \int_{-\infty}^{\infty} x^r\, e^{-x^2/2}\, dx $, for $3 \leq r \leq 10$).

$\bullet$ Derive experimentally (yet rigorously) {\bf explicit} expressions for many mixed moments with the more famous (and much more user-friendly) `measure if disarray', Netto's
{\it number of inversions}, and partially prove (by the method of moments) that these two permutation statistics are {\it jointly normal}, but,  not surprisingly, {\bf not independently so}
(unlike the pair (Number of Inversions, Major Index), see [BZ]). In fact we will prove that the asymptotic correlation between these two permutation statistics is $\frac{3}{\sqrt{10}}$,
and partially prove (by the method of moments) that the limiting distribution is bi-normal with that correlation, in other words, that the bi-variate pdf is
$$
\frac{1}{2 \,\sqrt{10}, \pi} \, \cdot \,  e^{
-\frac{x^2}{2} - \frac{y^2}{2} + \frac{3}{\sqrt{10}} \,xy
} \quad .
$$

{\bf A Quick Reminder about Weight-Enumerators and Probability Generating Functions of Discrete Random Variables Defined on Finite Sets, and how to Extract moments from them}

Let $A$ be a finite set and $f(a)$ some numerical attribute, called {\it random variable} (a very bad name!). We are interested in the average, variance, and higher moments of $f(a)$ over the
population $A$, assuming that we draw an element of $A$ {\it uniformly at random}. 
The {\bf weight-enumerator} according to $f(a)$ is the polynomial 
$$
\sum_{a \in A} q^{f(a)} \quad,
$$
and to get the {\bf probability generating function}, we simply divide by the number of elements of $A$. Calling that polynomial $F_A(q)$, we have:
$$
F_A(q) \, := \, \frac{1}{|A|} \, \sum_{a \in A} q^{f(a)} \quad .
$$
It is called the {\it probability generating function} since the coefficient of $q^i$ is the probability that the value of $f(a)$, where $a$ is drawn uniformly-at-random from $A$, happens to be $i$.

In terms of $F_A(q)$ we can easily get the {\bf expectation}, $\mu:=E[f]$, {\bf variance}, $m_2(f):=E[(f-\mu)^2]$,  and, more generally,
for $r>2$,  the so-called $r^{th}$ moment (about the mean), $m_r(f):=E[(f-\mu)^r]$, as follows.
$$
\mu=F_A'(1) \quad,
$$
and
$$
m_r(A) = \left ( q{\frac{d}{dq}} \right )^r \frac{F_A(q)}{q^\mu} {\Bigl \vert}_{q=1} \quad .
$$
Often we have an {\it infinite family} of natural sets, parameterized by positive integer $n$ (in this paper, the set of permutations on $\{1, \dots ,n\}$) and then
we are interested in {\bf explicit} expressions, in $n$, for the expectation, variance, and as many higher moments that we can get. See [Z], [BZ], and [CJZ] for nice examples.

{\bf How to get Explicit expressions to as many as possible moments of Spearman's Footrule?}

For the much more user-friendly permutation statistic $inv$, defined by
$$
inv(\pi):= \sum_{1 \leq i < j \leq n} [\pi_i > \pi_j] \quad,
$$
where $[statement]$ is $1$ or $0$  according to whether it is true or false, respectively, we have a beautiful, {\bf closed-form} expression for the  weight-enumerator, that goes back
(at least) to Netto:
$$
\sum_{\pi \in S_n} \, q^{inv(\pi)} \, = \, 1 \cdot (1+q) \cdots (1+q+ \dots + q^{n-1}) \, = \, \prod_{i=1}^{n} \frac{1-q^i}{1-q} \quad .
$$

This was treated (in the more general context of words) in [CJZ], but later we found out, thanks to Persi Diaconis, (see the erratum  to [CJZ]) that this was `old hat'. Nevertheless the
approach in [CJZ], of using symbolic computation to find explicit expressions for many moments, was novel, and gave a new, alternative,  elementary proof of this classical theorem.

We don't have this luxury for Spearman's footrule, since the weight-enumerator does not seem to have a nice  formula. 
Nevertheless, it is easy to see that for any $r$, the moment (and hence moment about the mean), is {\bf some} polynomial, and it
is also easy to bound its degree. So if we can compute many terms of the sequence of generating functions, let's call it $N_n(q)$
$$
N_n(q) \, := \, \sum_{\pi \in S_n} \, q^{na(\pi)} \quad,
$$
we can use the above process to get the moments numerically for sufficiently many $n$,  and finally fit  them with the appropriate polynomial.

If we had unlimited computer time and space, we can go as far as we wish by {\it pure brute force}. For example for finding the covariance (and hence correlation) of the
pair of permutation statistics $(inv,maj)$ (where $maj$ is the so-called {\it major index}), it worked very well. See [E], [K1], and [K2], Ex. 5.1.1.26).

Indeed, for any given $n$, there are only finitely many permutations of length $n$, so we can
find  $N_n(q)$ directly for any $n$. Alas, this is not very feasible starting at $n=11$. So we need to be more clever.

Let's write the permutation $\pi$ in {\it two-line-notation}.
$$
\left ( \matrix{ 1    &   2 & 3 & \dots & n \cr
                \pi_1 &   \pi_2 & \pi_3 & \dots & \pi_n}
\right ) \quad .
$$

Obviously
$$
na \left ( \matrix{ 1    &   2 & 3 & \dots & n \cr
                \pi_1 &   \pi_2 & \pi_3 & \dots & \pi_n}
\right ) \quad 
\, = \,
|\pi_1-1|+
na \left ( \matrix{   2 & 3 & \dots & n \cr
                   \pi_2 & \pi_3 & \dots & \pi_n}
\right ) \quad .
$$

So, in order to use {\it dynamical programming}, we are {\bf forced} to consider a more general creature
$$
N(S_1,S_2)(q):=\sum_{\pi \in FUN(S_1,S_2)}  q^{na(\pi)} \quad,
$$
where $S_1$ is a set of {\bf consecutive} integers, and $S_2$ is an {\it arbitrary} set of integers of the same cardinality, and $FUN(S1,S2)$ is the set of all one-to-one functions from $S_1$ to $S_2$.

The {\bf dynamical programming} natural recurrence is
$$
N\,(\{r,r+1,r+2, ..., r+k-1\} \, , \, \{a_1,\dots, a_k\})=
$$
$$
\sum_{i=1}^{k} 
q^{|r-a_i|} \cdot N( \{r+1,r+2, ..., r+k-1\},\{a_1,\dots, a_k\} \backslash \{a_i\}) \quad .
$$
This is implemented in procedure {\tt NPg} in the Maple package {\tt Noga12.txt}. Of course, {\it at the end of the day},  we only care about 
$$
N( \{1,2,3, ..., n\} , \{1,2,3, ..., n\} ) \quad,
$$
alias $N_n(q)$. This is procedure {\tt NPc(n,q)}. Procedure   {\tt NP(n,q)} does the same thing, straight from the definition, and of course should only be used for small $n$ (up to $n=9$) in order to
check the correctness of the `clever' procedure {\tt NPc}. Note that it still uses exponential time and memory, but $2^n$ is much smaller than $n!$, and we were able to compute
the first $25$ terms of $N_n(q)$. See the output file

{\tt https://sites.math.rutgers.edu/\~{}zeilberg/tokhniot/oNoga12c.txt} \quad .

It is easy to see that the expectation, variance, and higher moments (about the mean) , are {\it polynomials} in $n$. 
In fact that $r$-th moment about the mean has degree $\lfloor \frac{3r}{2} \rfloor$ in the variable $n$.
Diaconis and Graham used
human ingenuity to prove that
$$
\mu = \frac{(n-1)(n+1)}{3} 
$$
$$
\sigma^2 = \frac{(n+1)\,(2n^2+7)}{45} \quad ,
$$
but did not derive higher moments. Using our experimental-yet-rigorous method, we derived explicit expressions for all the moments up to the $10^{th}$.
$$
m_3 \,  = \, -\frac{2 \left(n +2\right) \left(n +1\right) \left(2 n^{2}+31\right)}{945}
$$

$$
m_4 \,  = \, \frac{\left(n +1\right) \left(28 n^{5}+180 n^{3}+160 n^{2}+887 n +1265\right)}{4725}
$$

$$
m_5 \,  = \, -\frac{4 \left(n +2\right) \left(n +1\right) \left(44 n^{5}-10 n^{4}+788 n^{3}+86 n^{2}+3587 n +8555\right)}{93555}
$$

$$
m_6 \,  = \, \frac{\left(n +1\right)}{127702575} \cdot
$$
$$
(168168 n^{8}-145288 n^{7}+1800148 n^{6}+2180892 n^{5}+18508182 n^{4}+32547228 n^{3}
$$
$$
+112117257 n^{2}+385870348 n +368963105 ) \quad .
$$

For explicit expressions for $m_7,m_8,m_9$ and $m_{10}$ see the output file

{\tt https://sites.math.rutgers.edu/\~{}zeilberg/tokhniot/oNoga12a.txt} \quad .

From these expressions, Maple  (and even you!) can compute the limits of the {\bf scaled moments about the mean}, $\frac{m_r}{m_2^{r/2}}$, and
see that for $3 \leq r \leq 10$ they coincide with those of the standard normal distribution $\N(0,1)$ (namely $0,3,0,15,0,105,0,945$),  giving ample evidence that
Spearman's footrule is {\it asymptotically} normal, and rigorously proving it up to the $10^{th}$ moment. Of course, as mentioned above this  was already proved, fully,  by Persi Diaconis in his
joint paper with Ron Graham [DG], but he used `fancy stuff'.

{\bf How does Spearman's Footrule ``interact'' with Netto's number of inversions}

The dynamical programming recurrence that enabled us to compute the first $25$ terms of the sequence of weight-enumerators $N_n(q)$ can be easily modified to
compute many terms of the {\bf bi-variate} polynomial, let's call it $S_n(p,q)$

$$
S_n(p,q) \, := \, \sum_{\pi \in S_n} \, p^{inv( \pi )} q^{na(\pi)} \quad,
$$
see procedure  {\tt NPcJ(n,p,q)} in our Maple package. This enabled us to compute the first $22$ terms,  that can be viewed in the following output file

{\tt https://sites.math.rutgers.edu/\~{}zeilberg/tokhniot/oNoga12d.txt } \quad .

From these we immediately derived, by the same empirical method, the covariance
$$
\frac{(n+1) \cdot (n^2+1)}{30} \quad,
$$
from which followed that the {\bf asymptotic correlation} is  $\frac{3}{\sqrt{10}}$, and  we derived the {\bf mixed moments}, $m_{r,s}:=E[(inv-\mu_{inv})^r (na- \mu_{na})^s]$ for all $r+s \leq 8$.

To see them for $1 \leq r,s \leq 4$ see the output file

{\tt https://sites.math.rutgers.edu/\~{}zeilberg/tokhniot/oNoga12b.txt} \quad,

from which it followed that the {\bf scaled mixed moments} (at least for $r+s \leq 8$) tend to those of the bi-variate normal distribution with
the above-mentioned correlation, namely the one whose joint pdf (probability density function) is
$$
\frac{1}{2 \,\sqrt{10}\, \pi} \, \cdot \,  e^{-\frac{x^2}{2} - \frac{y^2}{2} + \frac{3}{\sqrt{10}} \,xy
} \quad .
$$

This gives ample evidence that it is true for {\it all} mixed moments, and that, in turn,  would conclusively prove (using the method of moments) that the pair $(inv,na)$ is
jointly asymptotically normal with correlation $\frac{3}{\sqrt{10}}$. 
We are wondering whether Persi Diaconis, or any of his disciples, can
prove this fact with their heavy machinery.

{\bf Postscript} (written Oct. 4, 2021). After the first version was posted, and one of (DZ) gave a talk about it at the Rutgers Experimental Mathematics Zoom seminar, 
Stoyan Dimitrov pointed our attention to a very interesting paper 

``The generating function for total displacement'' by T. Kyle Petersen and Mathieu Guay-Paquet (Elec. J. Combinatorics {\bf 21} (2014))

{\tt https://arxiv.org/pdf/1404.4674.pdf}  \quad,

that enables a much faster computation of the weight-enumerators $N_n(q)$. Their approach enables a fast computation of what we called $N_n(q)$ 
for $1 \leq n \leq 50$, that can be gotten by typing 

{\tt SeqF(n,50)},

in the new version of our Maple package, where we added a new procedure, {\tt SeqF}, implementing their nice approach (using weighed Motzkin paths, and the implied `infinite' continued fraction).
Using this more extensive data set, we were able to find the first $19$ moments, see the new output file:

{\tt https://sites.math.rutgers.edu/\~{}zeilberg/tokhniot/oNoga12e.txt} \quad. 

We also learned, thanks to that paper, that Spearman's footrule is mentioned 
in Knuth's magnum-opus [K2], as exercise 5.1.1.28. Also Martin Rubey pointed out that nowadays one does not have to bother Persi Diaconis, all one has to do is
consult the wonderful website

{\tt https://www.findstat.org/} \quad . See his following insightful message:

{\tt https://sites.math.rutgers.edu/\~{}zeilberg/mamarim/mamarimhtml/MartinRubeyComments.pdf}

Of related interest is the following interesting paper: `Moments of permutation statistics and central limit theorems', by
Stoyan Dimitrov and Niraj Khare, available from

{\tt https://arxiv.org/abs/2109.09183} \quad .

{\bf PostPostscript} (written Oct. 7, 2021): We now read the above mentioned lovely article by
T. Kyle Petersen and Mathieu Guay-Paquet more carefully, and realized that using the continued fraction, while very elegant, is not
the most efficient way to generated many terms of $\{N_n(q)\}$. What we called before {\tt SeqF(q,N)} (to generate the first {\tt N} terms
of that sequence) has been now renamed {\tt SeqFcf(q,N)}, and the new procedure, {\tt SeqF(q,N)}, uses the fact that $N_n(q)=F(n,0)$
where $F(a,b)$ is the weight enumerators of truncated Motzkin paths (with their weight), that terminate at the point $(a,b)$.

It is easy to see (once Petersen and Guay-Paquet introduced their beautiful approach) that, $F(a,b)$ satisfies the {\bf dynamical programming} recurrence (but
no longer exponential time and memory, but only quadratic time and memory):
$$
F(a,b) \, = \, q^{2b} \left (\, (b+1)  F(a-1,b+1)\, +\, b F(a-1,b-1)+ (2b+1)  F(a-1,b) \, \right ) \quad,
$$
subject to the initial condition $F(0,0)=1$ and the boundary condition $F(a,b)=0$ if $b<0$.

{\bf PostPostPostscript} (written Oct. 14, 2021): T. Kyle Petersen kindly told us (private communication) how to $q$-analogize (but we use $p$ for inversions here) the above recurrence
to efficiently compute what we call above $S_n(n,p,q)$. Definte $G(a,b)$ to be the polynomials in $(p,q)$ satisfying the recurrence
$$
G(a,b) \, = \, (pq^2)^{b} \left (\, [b+1]  G(a-1,b+1)\, +\, [b] G(a-1,b-1)+ ([b]+[b+1])  F(a-1,b) \, \right ) \quad,
$$
subject to the initial condition $G(0,0)=1$ and the boundary condition $G(a,b)=0$ if $b<0$.

Here (as usual) $[b]=1+p+ \dots + p^{b-1}$. Then $S_n(p,q)=G(n,0)$.

This is now implemented in procedure {\tt SeqFg}. For example to get the first $22$ terms, type {\tt SeqFg(p,q,22);}.

{\bf References}

[BZ] Andrew Baxter and Doron Zeilberger,
{\it The Number of Inversions and the Major Index of Permutations are Asymptotically Joint-Independently-Normal (Second Edition)},
The Personal Journal of Shalosh B. Ekhad and Doron Zeilberger.  Feb. 4, 2011. \hfill\break
{\tt https://sites.math.rutgers.edu/\~{}zeilberg/mamarim/mamarimhtml/invmaj.html} \quad. Also published in
{\tt https://arxiv.org/abs/1004.1160}.

[CJZ] E. Rodney Canfield, Svante Janson, and Doron Zeilberger, {\it The Mahonian probability distribution  on words IS asymptotically normal},
Adv. Appl. Math. {\bf 46} (2011), 109-124. \hfill\break
{\tt https://sites.math.rutgers.edu/\~{}zeilberg/mamarim/mamarimhtml/mahon.html} \quad .\hfill\break
Erratum: \hfill\break
{\tt https://sites.math.rutgers.edu/\~{}zeilberg/mamarim/mamarimPDF/mahonerratum.pdf} \quad .

[DG] Persi Diaconis and Ron Graham, {\it Spearman's Footrule as a Measure of Disarray},Journal of the Royal Statistical Society: Series B (Methodological) {\bf 39} (1977), 262-268.

[E] Shalosh B. Ekhad, 
{\it The joy of brute force: the covariance of the major index and the number of inversions}, Personal Journal of S. B. Ekhad and D. Zeilberger, , ca. 1995.\hfill\break
{\tt https://sites.math.rutgers.edu/\~{}zeilberg/mamarim/mamarimhtml/brute.html}

[K1] Donald E. Knuth, {\it Letter to Doron Zeilberger}, May 16, 1995. \hfill\break
{\tt https://sites.math.rutgers.edu/\~{}zeilberg/mamarim/mamarimhtml/dkMay95.pdf} \quad .

[K2] Donald E. Knuth, {\it The Art of Computer Programming, Volume 3, Sorting and Searching, Second Edition}, Addison-Wesley, 1998.

[W] Jessica Wynne, {\it ``Do Not Erase, Mathematicians and their chalkboard''}, Princeton University Press, 2021.

[Z] Doron Zeilberger, {\it HISTABRUT: A Maple Package for Symbol-Crunching in Probability theory}, The Personal Journal of Shalosh B. Ekhad and Doron Zeilberger, Aug. 25, 2010. \hfill\break
{\tt https://sites.math.rutgers.edu/\~{}zeilberg/mamarim/mamarimhtml/histabrut.html} \quad. Also published in:
{\tt https://arxiv.org/abs/1009.2984} \quad .

\bigskip
\hrule
\bigskip
Shalosh B. Ekhad and Doron Zeilberger, Department of Mathematics, Rutgers University (New Brunswick), Hill Center-Busch Campus, 110 Frelinghuysen
Rd., Piscataway, NJ 08854-8019, USA. \hfill\break
Email: {\tt [ShaloshBEkhad, DoronZeil] at gmail dot com}   \quad .

{\bf Exclusively published in the Personal Journal of Shalosh B. Ekhad and Doron Zeilberger and arxiv.org}

First Written: {\bf Sept. 29, 2021}. Second version: (thanks to Stoyan Dimitrov and Martin Rubey): {\bf Oct. 4, 2021}.
Third version: {\bf Oct. 7, 2021} \quad . This version (thanks to T. Kyle Petersen): {\bf Oct. 14, 2021} \quad .

\end